\documentclass[journal]{IEEEtran}

\ifCLASSINFOpdf

\else

\fi

\usepackage{amsmath,amssymb,amsthm}
\usepackage{graphicx}
\usepackage{cite}
\usepackage{color}
\usepackage{tabularx}
\usepackage{hyperref}
\usepackage{tikz}
\usetikzlibrary{arrows,decorations.pathmorphing,backgrounds,fit,positioning,shapes.symbols,chains}

\makeatletter
\g@addto@macro{\UrlBreaks}{\UrlOrds}
\makeatother

\begin{document}

\title{A Robust Approach to Chance Constrained Optimal Power Flow with Renewable Generation}

\author{Miles Lubin, Yury Dvorkin, \textit{Student Member, IEEE,} Scott Backhaus

\thanks{Manuscript received April 23, 2015; revised August 7, 2015.
}
\thanks{Miles Lubin is with the Operation Research Center, Massachusetts Institute of Technology, Cambridge, MA, 02142 USA (e-mail: \href{mailto:mlubin@mit.edu}{mlubin@mit.edu}).
}
\thanks{Yury Dvorkin is with the Department of Electrical Engineering, University of Washington, Seattle, WA, 98195 USA (e-mail: \href{mailto:dvorkin@uw.edu}{dvorkin@uw.edu}).
}

\thanks{Scott Backhaus is with the Center for Nonlinear Studies, Los Alamos National Laboratory, Los Alamos, NM, 87545 USA (e-mail: \href{mailto:backhaus@lanl.gov}{backhaus@lanl.gov}). 
}
}

\markboth{Submitted to IEEE Transactions on Power Systems}%
{Shell \MakeLowercase{\textit{et al.}}: Bare Demo of IEEEtran.cls for Journals}

\maketitle

\begin{abstract}
Optimal Power Flow (OPF) dispatches controllable generation at minimum cost subject to operational constraints on generation and transmission assets. The uncertainty and variability of intermittent renewable generation is challenging current deterministic OPF approaches. Recent formulations of OPF use chance constraints to limit the risk from renewable generation uncertainty, however, these new approaches typically assume the probability distributions which characterize the uncertainty and variability are known exactly. We formulate a Robust Chance Constrained (RCC) OPF that accounts for uncertainty in the parameters of these probability distributions by allowing them to be within an uncertainty set. The RCC OPF is solved using a cutting-plane algorithm that scales to large power systems. We demonstrate the RRC OPF on a modified model of the Bonneville Power Administration network, which includes 2209 buses and 176 controllable generators. Deterministic, chance constrained (CC), and RCC OPF formulations are compared using several metrics including cost of generation, area control error, ramping of controllable generators, and occurrence of transmission line overloads as well as the respective computational performance.
\end{abstract}


\IEEEpeerreviewmaketitle

\section{Introduction}
\subsection{Motivation}
The continually growing penetration of intermittent renewable energy resources, e.g., wind and solar photovoltaic, is revealing a number of drawbacks in existing power system operational procedures that may limit the integration of these new resources. Wind generation is an intermittent and not fully dispatchable generation technology that imposes challenges to least-cost, risk-averse management of generation and transmission assets. One approach to address these challenges is strategic investments  in more transmission and controllable generation capacity to enhance the system flexibility \cite{ma_2013,ruth_2015}. These investments are costly and subject to a variety of regulatory and policy limitations. On the other hand, improving operating protocols may create additional flexibility in the existing system by replacing ad hoc deterministic policies for limiting system risk from intermittent generation with probabilistic formulations that account for intermittency in a principled manner. 

Historically, these deterministic policies were designed to account for less challenging deviations of load from its forecast value. They have performed very well for their design conditions where fluctuations in load are a small fraction of the total load, however, they are not expected to cost effectively manage risk when net-load fluctuations are large. Recent regulatory initiatives, such as Federal Energy Regulatory Commission (FERC) Orders 764 \cite{ferc_764} and 890 \cite{ferc_890} have identified the need for a new generation of operating protocols and decision-making tools for the successful integration of renewable generation. 

In this manuscript, we implement a distributionally robust chance constraint (RCC) optimal power flow (OPF) model and compare it with the deterministic OPF and chance constrained (CC) OPF models. The deterministic OPF is a typical short-term decision-making tool used by a number of utilities and its implementation in this work aims to give a reasonable benchmark for comparison. The chance constrained (CC) OPF limits the probability of violating transmission or generation constraints using a statistical model of \textcolor{black}{wind deviations from forecast values, which are parametrized with a zero-mean Gaussian distribution \cite{yli_2014, letter}}. As compared to the CC OPF, the RCC OPF is a generalization that allows for uncertainty in the mean and variance of the wind forecast error. To demonstrate the effect of using probabilistic methods in an OPF, we compare these different OPF formulations in the setting of vertically-integrated grid operations, specifically, on a modification of the Bonneville Power Administration (BPA) system. 

\subsection{\textcolor{black}{Literature Review}}

Traditional deterministic OPF models \cite{bacher_1988} dispatch controllable generation using the central (most likely) wind forecast, i.e., they do not endogenously account for the variability and uncertainty of wind generation \cite{ela_2012}. Exogenously calculated reserve margins and heuristic policies \cite{lew_2010,doherty_2005} are often used to enforce additional security requirements to ensure system reliability. However, these heuristic approaches are limited in their ability to produce cost-efficient solutions \cite{ruiz_2009}, \cite{morales_2009}. 

Recently, a number of transmission-constrained OPF and unit commitment (UC) models based on stochastic programming \textcolor{black}{\cite{yli_2014,oren_2011,wang2008,ryan_2013}}, interval programming \cite{yang_2011},\cite{lei_2012}, chance constrained optimization \cite{ccopf-sirev} and robust optimization \textcolor{black}{\cite{rjabr_2013b, jabr_2014, bertsimas_2013}} have been proposed for endogenous risk-averse decision making. A common drawback of stochastic programming and chance constrained optimization is a requirement for accurate statistical models of wind generation uncertainty and variability. In practice, wind generation is modeled using the uni- or multivariate Gaussian distribution, which unavoidably results in solution inaccuracy \cite{pinson_2010}. In addition, stochastic programming typically requires the generation of  a relatively large number of scenarios leading to impractically large computing times, even for relatively large values of the duality gap \cite{papavasiliou_2014}. If the number of scenarios is reduced by means of scenario reduction techniques \cite{dupacova_2000}, the monetary benefits attained with stochastic programming may reduce accordingly \cite{lei_2012},\cite{papavasiliou_2014}.

In contrast, interval programming and robust optimization models allow wind generation fluctuations within a certain range around a central forecast. By disregarding the likelihood of individual scenarios within this range, these methods result in an overly conservative solution as compared to stochastic programming \cite{guan_2013}. Although the conservatism of robust optimization can be reduced by adjusting the budget of uncertainty \cite{bertsimas_2004} or using  dynamic uncertainty sets \cite{lorca_2014}, there is no systematic methodology to choose the value of budget of uncertainty \textit{a priori}.

Here, we leverage recent work by Bienstock et al~\cite{ccopf-sirev} that formulated and implemented a CC OPF using a cutting-plane approach and demonstrated its scalability to large systems. Based on the CC OPF,  Bienstock et al~\cite{ccopf-sirev} also envisioned an RCC OPF that allowed for uncertainties in the mean and variance of wind forecast errors, but did not implement or test the RCC OPF. In this manuscript, we build upon the RCC OPF formulation from \cite{ccopf-sirev} to reduce the inaccuracy of assuming Gaussian distributions for wind forecast errors. The RCC OPF is implemented in JuMPChance \cite{jumpchance}, an open source optimization package developed to model ordinary and distributionally robust chance constraints in OPF and other optimization settings. We use a modification of the BPA system to demonstrate the performance of the RCC OPF implementation on large-scale interconnections. Our case study shows that robustification of chance constraints can result in cost savings and is a more effective decision-making tool to mitigate real-time active power imbalances, extreme ramping of generators and transmission overloads. 

The rest of this manuscript is organized as follows. Section~\ref{sec:RCCOPF} reviews a deterministic OPF, restates (with a few modifications) the CC OPF of \cite{ccopf-sirev}, and demonstrates how  robustified chance constraints and are incorporated into the proposed RCC OPF model. Section~\ref{sec:case_study} describes the BPA test system and compares the performance of the deterministic OPF, CC OPF, and RCC OPF models on this test system. Finally, \textcolor{black}{Section}~\ref{sec:conclusions} presents conclusions and possible directions for future work.

\section{RCC OPF Formulation}\label{sec:RCCOPF}

\textcolor{black}{Section~\ref{sec:CCOPF_notations} defines notations used throughout Section~\ref{sec:RCCOPF}.} In \textcolor{black}{Section}~\ref{sec:CCOPF_formulation}, we review the CC OPF formulation from \cite{ccopf-sirev} and restate it with a few modifications.  Section~\ref{sec:robustmodel} discusses a model for distributionally robust chance constraints which are incorporated into an RCC OPF formulation in Section~\ref{sec:RCC-OPF_formulation}.

\subsection{Notations}
\label{sec:CCOPF_notations}
\noindent$\mathcal{B}$ -- set of buses\\
$r$ -- index of reference bus ($\in \mathcal{B}$)\\
$\mathcal{L}$ -- set of lines\\
$\mathcal{G}$ -- set of controllable generators\\
$G_b$ -- subset $(\subset \mathcal{G})$ of generators located at bus $b$\\
$\mathcal{W}$ -- subset of buses with wind farms\\
$\beta_{mn}$ -- susceptance of line $(m,n)$\\
$f_{mn}$ -- real power flow over line $(m,n)$, MW\\
$\theta_b$ -- phase angle at bus $i$\\
$p_i$ -- output of controllable generator $i$, MW\\
$d_b$ -- demand at bus $b$, MW\\
$p_i^{min}$ -- minimum output of generator $i$, MW\\
$p_i^{max}$ -- maximum output of generator $i$, MW\\
$f_{mn}^{max}$ -- capacity of line $(m,n)$, MW\\
$RU_i$ -- max ramp-up of generator $i$ in the OPF period, MW/h\\
$RD_i$ -- max ramp-down of generator $i$ in the OPF period, MW/h\\
$w^f_b$ -- forecast output of wind farm at bus $b$, MW\\
$\omega_b(t)$ -- actual deviation from forecast $w_b^f$ at time $t$, MW\\
$c_{i1}$ -- linear coefficient of cost for generator $i$, \$/MW \\
$c_{i2}$ -- quadratic coefficient of cost for generator $i$, \$/MW$^2$\\
\textcolor{black}{$\delta$ -- nonphysical auxiliary variable}\\
\textcolor{black}{$B$ -- bus admittance matrix}\\
\textcolor{black}{$\pi_b$ -- $b$th row of the inverse of the admittance matrix, after excluding the row and column corresponding to the reference bus}\\

\noindent In the rest of this manuscript, bold symbols denote random variables. In particular, $\boldsymbol\omega_b$ models deviations $\omega_b(t)$ within the OPF period, which drive random fluctuations in controllable generator injections $\boldsymbol p_i$, bus phases $\boldsymbol \theta_b$, and power line flows $\boldsymbol f_{mn}$ (described below). We denote the total deviation from the forecast as $\boldsymbol\Omega = \sum_{b \in \mathcal{W}} \boldsymbol\omega_b$. In the CC OPF, the deviations  $\boldsymbol\omega_b$ are assumed \textit{independent and normally distributed with zero mean and known variance $\sigma_b^2$}. In the RCC OPF formulation discussed in \textcolor{black}{Section ~\ref{sec:RCC-OPF_formulation}}, this assumption is relaxed by introducing symmetric intervals $[-\bar \mu_b,\bar\mu_b]$ and $[\sigma_b^2 - \bar\sigma_b^2,\sigma_b^2 + \bar\sigma_b^2]$ for the mean and variance of $\boldsymbol\omega_b$.

\subsection{CC OPF Formulation}\label{sec:CCOPF_formulation}
The CC OPF formulation is derived from the following single-stage deterministic OPF:
\begin{equation}
\min_{p,\theta} \sum_{i \in \mathcal{G}} (c_{i2}p_i^2 + c_{i1}p_i)\label{eq:Det_OPF}
\end{equation}
\text{subject to} 
\begin{gather}
\sum_{n\in \mathcal{B}} B_{bn}\theta_n = \sum_{i \in G_b} p_i + w^f_b - d_b, \quad \forall b \in \mathcal{B}, \\
p_{i}^{min} \leq p_i \leq p_i^{max}, \quad \forall i \in \mathcal{G},\\
\label{eq:flowdef} f_{mn} = \beta_{mn}(\theta_m - \theta_n), \quad \forall \{m,n\} \in \mathcal{L}, \\
\label{eq:Det_OPF_end} -f_{mn}^{max} \leq f _{mn}  \leq f_{mn}^{max}, \quad \forall \{m,n\} \in \mathcal{L}, 
\end{gather}
where $B$ is the $|\mathcal{B}| \times |\mathcal{B}|$ bus admittance matrix defined by:
\begin{equation}
B_{mn} = \left\{
\begin{array}{cc} -\beta_{mn}, & \{m,n\} \in \mathcal{L}\\
\sum_{k:\{k,n\}\in\mathcal{E}} \beta_{kn}, &  m = n\\
0, & \text{otherwise.}\end{array}\right.
\end{equation}

In the deterministic OPF formulation in \eqref{eq:Det_OPF}-\eqref{eq:Det_OPF_end}, the controllable generation set points $p_i$ are optimized to minimizing the total cost of generation for the forecast wind generation $\omega^f$  and demand $d_{b}$ subject to operating constraints on generators and transmission lines. In the presence of deviations $\boldsymbol\omega_b$ from $\omega_b^f$, we model the proportional response of generators:
\begin{equation}
\label{eq:propresponse} 
\boldsymbol p_i = p_i - \alpha_i \boldsymbol\Omega. 
\end{equation}
Here, $\alpha_i\geq 0$ is the participation factor for controllable generator $i$. When $\sum_i \alpha_i=1$, the response rule \eqref{eq:propresponse} guarantees that generation and load remain balanced, but does not limit the magnitude of the response of the generators or the resulting flow on the power lines. 

As shown in \cite{ccopf-sirev}, the deterministic OPF in \eqref{eq:Det_OPF} can be reformulated as a CC OPF by introducing probabilistic constraints on $\boldsymbol f_{mn}$ and $\boldsymbol p_i$ and modelling the participation factors $\alpha_i$ as decision variables. \textcolor{black}{We present the CC OPF formulation as follows:}
\begin{align}
& \min &  \sum_{i \in \mathcal{G}}\left(c_{i2}(p_i^2 + \operatorname{var}(\boldsymbol\Omega)\alpha_i^2) + c_{i1}p_i\right) \label{eq:CC-OPF_objective} 
\end{align}
\text{subject to:} 
\begin{gather}
 \sum_{n\in \mathcal{B}} B_{bn}\theta_n = \sum_{i \in G_b} p_i + w^f_b - d_b, \quad \forall b \in \mathcal{B}\label{eq:CC-OPF_PF} \\
 p_{i}^{min} \leq p_i \leq p_i^{max}, \quad \forall i \in \mathcal{G} \\
 f_{mn} = \beta_{mn}(\theta_m - \theta_n), \quad \forall \{m,n\} \in \mathcal{L} \\
 |f_{mn}| \leq f_{mn}^{max}, \quad \forall \{m,n\} \in \mathcal{L}\label{eq:CC-OPF_linelimit} \\
 P(p_i - \boldsymbol\Omega \alpha_i \le p_i^{max}) \ge 1- \epsilon_i, \quad \forall i \in \mathcal{G}\label{eq:CC-OPF_pmax} \\
 P(p_i - \boldsymbol\Omega \alpha_i \ge p_i^{min}) \ge 1- \epsilon_i, \quad \forall i \in \mathcal{G}\label{eq:CC-OPF_pmin} \\
 P(-\boldsymbol\Omega\alpha_i \le RU_i) \ge 1-\epsilon_i, \quad \forall i \in \mathcal{G}\label{eq:CC-OPF_RU} \\
 P(\boldsymbol\Omega\alpha_i  \le RD_i) \ge 1- \epsilon_i, \quad \forall i \in \mathcal{G}\label{eq:CC-OPF_RD} \\
P\bigg(f_{mn} + \beta_{mn}\boldsymbol\Omega(\delta_n - \delta_m) + \beta_{mn}\boldsymbol\omega^T(\pi_m-\pi_n) \nonumber \\ \le f_{mn}^{max} \bigg)  \ge 1- \epsilon_{mn}, \quad \forall \{m,n\} \in \mathcal{L}\label{eq:CC-OPF_fplus} \\
P\bigg( f_{mn} + \beta_{mn}\boldsymbol\Omega(\delta_n - \delta_m) +\beta_{mn}\boldsymbol\omega^T(\pi_m-\pi_n) \nonumber \\ \ge -f_{mn}^{max} \bigg) \ge 1- \epsilon_{mn}, \quad \forall \{m,n\} \in \mathcal{L}\label{eq:CC-OPF_fminus} \\
\sum_{i\in\mathcal{G}} \alpha_i = 1, \quad \alpha \geq 0 \\
\sum_{i\in G_r}\alpha_i = \delta_r = \theta_r = 0.\label{eq:CCOPF_last}  \\
\sum_{\substack{n\in \mathcal{B}\\n\neq r}} B_{bn}\delta_n = \textcolor{black}{\sum_{i \in G_b} \alpha_i}, \quad \forall b \in \mathcal{B}\setminus \{r\}\label{eq:CC-OPF_delta} 
\end{gather}
Here the decision variables are $p, \theta, \delta, \alpha$ and $f$.  For the quadratic cost of production (assumed to be convex), \cite{ccopf-sirev} shows that
\begin{equation}\label{eq:ccobjective}
\mathbb{E}[\sum_{i \in \mathcal{G}} c_{i2}\boldsymbol p_i^2 + c_{i1}\boldsymbol p_i] = \sum_{i \in \mathcal{G}}\left(c_{i1}(p_i^2 + \operatorname{var}(\boldsymbol\Omega)\alpha_i^2) + c_{i2}p_i\right).
\end{equation}

Therefore, in (\ref{eq:CC-OPF_objective}), the CC OPF seeks to minimize the convex quadratic expected cost of production.  Under these assumptions, the CC OPF is tractable and representable using second-order cone programming (SOCP).

Constraints (\ref{eq:CC-OPF_PF})-(\ref{eq:CC-OPF_linelimit}) are deterministic and enforce power flow feasibility, generator limits, and power line flow limits for $\omega^f$, similar to those in the deterministic OPF of (\ref{eq:Det_OPF}). In chance constraints (\ref{eq:CC-OPF_pmax})-(\ref{eq:CC-OPF_pmin}), the controllable generator outputs are now random and given by Eq.~\eqref{eq:propresponse}. As in \cite{ccopf-sirev}, these chance constraints bounded by $\epsilon_i$, the probability of the fluctuating generator outputs exceeding their upper limits $p_i^{max}$ or lower limits $p_i^{min}$.

\textcolor{black}{Relative to the CC OPF formulation of \cite{ccopf-sirev}, we add only the constraints}~\eqref{eq:CC-OPF_RU} and \eqref{eq:CC-OPF_RD} that limit the probability of the real-time response $\alpha_i \boldsymbol\Omega$ to wind deviations from the forecast from exceeding $RU_i$ for positive changes and $RD_i$ for negative changes.
Here, $RU_i$ and $RD_i$ are the continuous ramping constraints on the generators over the OPF time step.

Chance constraints \eqref{eq:CC-OPF_pmax}-\eqref{eq:CC-OPF_RD} on the controllable generator injections are expressed explicitly in terms of $\alpha$ and the random wind fluctuations $\boldsymbol \omega$.  Chance constraints \eqref{eq:CC-OPF_fplus}-\eqref{eq:CC-OPF_fminus}, which bound the line flows $\boldsymbol f_{mn}$, are more subtle.  The flows $\boldsymbol f_{mn}$ also change with the fluctuating wind injections and the controllable generator response, however, the $\boldsymbol f_{mn}$ depend on the wind deviations in an implicit manner, i.e.
$
\boldsymbol f_{mn} = \beta_{mn}(\boldsymbol\theta_m - \boldsymbol\theta_n),
$
where
\begin{equation}\label{eq:boldthetasystem}
\sum_{n\in \mathcal{B}} B_{bn}\boldsymbol\theta_n = \sum_{i \in G_b} (p_i -\alpha_i\boldsymbol\Omega) + w^f_b - d_b + \boldsymbol\omega_b \quad\forall\, b \in \mathcal{B}.
\end{equation}
Bienstock et al.~\cite{ccopf-sirev} derive explicit equations for $\boldsymbol f_{mn}$ by observing that once a reference bus $r \in \mathcal{B}$ is chosen and $\boldsymbol\theta_r$ and $\alpha_r$ fixed to zero, the system of equations~\eqref{eq:boldthetasystem} is invertible and the adjusted phase angles $\boldsymbol\theta$ (and hence $\boldsymbol f$) can be expressed as a linear function of $\theta$, $\boldsymbol\Omega\delta$, and $\boldsymbol\omega$:
\begin{equation}\label{eq:boldtheta}
\boldsymbol\theta_b = \theta_b - \boldsymbol\Omega\delta_b + \pi_b^T\boldsymbol\omega
\end{equation}
where $\theta$ satisfies \eqref{eq:CC-OPF_PF}, $\delta$ satisfies \eqref{eq:CC-OPF_delta}, and
$\pi_b$ is the $b$th row of $\hat B^{-1}$ (oriented as a column vector), where $\hat B$ is the ($|\mathcal{B}| - 1 \times |\mathcal{B}| - 1$) submatrix of $B$ with the row and column corresponding to the reference bus removed.  The variable $\delta$ is introduced solely for computational convenience. 
The chance constraints \eqref{eq:CC-OPF_fplus}-\eqref{eq:CC-OPF_fminus} use~\eqref{eq:boldtheta} to express $\boldsymbol f_{mn}$ explicitly in terms of the random variables and decision variables.

While chance constraints like \eqref{eq:CC-OPF_pmax}-\eqref{eq:CC-OPF_fminus} are often nonconvex and difficult to treat in general~\cite{NemirovskiShapiro07}, under the assumption of normality, they are both convex and computationally tractable~\cite{PrekopaChapter}. In particular,
a chance constraint of the form
\begin{equation}\label{eq:prob}
\mathbb{P}(\boldsymbol\xi ^T x \leq b) \geq 1-\epsilon
\end{equation}
is equivalent to
\begin{equation}\label{eq:reform}
\mu^Tx + \Phi^{-1}(1-\epsilon)\sqrt{x^T\Sigma x} \leq b,
\end{equation}
when $\boldsymbol\xi \sim N(\mu,\Sigma)$ where $\Phi^{-1}$ is the inverse cumulative distribution function of the standard normal distribution.
In the following, we assume $\epsilon < \frac{1}{2}$ so that   $\Phi^{-1}(1-\epsilon) > 0$ and constraint \eqref{eq:reform} is convex in $(x,b)$. Note that in this model, we treat each chance constraint independently.
Although it would also be natural to pose a model which attempts to enforce that multiple linear constraints hold jointly with high probability, convexity in this case remains an open question, even under the assumption of normality~\cite{henrion_2005,NemirovskiShapiro07}.

Constraint~\eqref{eq:reform} is not only convex; it can be represented as a second-order cone (SOC) constraint handled by many off-the-shelf optimization packages like CPLEX~\cite{cplex} and Gurobi~\cite{gurobi}. Indeed, we see that~\eqref{eq:reform} is
satisfied iff \textcolor{black}{$\exists t$ such that}
\begin{gather}
\color{black}
t \ge ||\Sigma^{\frac{1}{2}} x||_2,\label{eq:slacksoc} \\
\color{black}
\mu^Tx + \Phi^{-1}(1-\epsilon)t \le b,
\end{gather}
\textcolor{black}{where~\eqref{eq:slacksoc} is a standard SOC constraint.}

While the representation of the CC OPF in~\eqref{eq:CC-OPF_objective}-\eqref{eq:CCOPF_last} 
with the reformulation of the chance constraints per \eqref{eq:reform} is quite useful, 
\cite{ccopf-sirev} observed that off-the-shelf solvers were not capable of solving large scale CC OPF instances. Instead, \cite{ccopf-sirev} implemented a specialized algorithm based on sequential \textcolor{black}{outer approximation of~\eqref{eq:slacksoc}. We note that in our implementation,} the formulation is provided as stated in~\eqref{eq:CC-OPF_objective}-\eqref{eq:CCOPF_last} to the modeling tool JuMPChance~\cite{jumpchance}, \textcolor{black}{which enables the user to select between solution via sequential outer approximation or via reformulation to SOCP.}

\subsection{Cutting-plane algorithm to solve distrbutionally robust chance constraints}\label{sec:robustmodel}

In the analytic reformulation of the chance constraints in \eqref{eq:reform}, the wind deviations $\boldsymbol \omega_b$ are assumed to be normally distributed with known (zero) mean and variance, i.e. $\boldsymbol\omega \sim N(\mu,\Sigma)$.
This approach is computationally tractable but has drawbacks. The assumption of normality is often an approximation, and even when valid, $\mu$ and $\Sigma$ are typically estimated from data and not known exactly. Often, we can only say with confidence that $(\mu,\Sigma)$ fall in some uncertainty set $U$. By reformulating the chance constraints of \textcolor{black}{Section}~\ref{sec:CCOPF_formulation} to so-called distributionally robust chance constraints, the constraint $\xi^Tx \leq b$ is required to hold with high probability \textit{under all possible distributions} within $U$, i.e.,
\begin{equation}\label{eq:robustprob}
\mathbb{P}_{\xi \sim N(\mu,\Sigma)}(\xi ^Tx \leq b) \geq 1-\epsilon,\quad \forall\, (\mu,\Sigma) \in U.
\end{equation}

For each $(\mu,\Sigma)\in U$, we have a single convex constraint of the form~\eqref{eq:prob}, therefore~\eqref{eq:robustprob} is a potentially infinite set of convex constraints and is convex itself. \textcolor{black}{ In this manuscript, as in~\cite{ccopf-sirev}, we consider uncertainty sets $U$ that can be partitioned into a product $U = U_\mu \times U_{\Sigma}$ where $(\mu,\Sigma) \in U$ iff $\mu \in U_\mu$ and $\Sigma \in U_\Sigma$.
Under these assumptions, Bienstock et al.~\cite{ccopf-sirev} observe that there is no known compact, deterministic reformulation of~\eqref{eq:robustprob} and instead describe a cutting-plane algorithm which we build upon and demonstrate is capable of handling large-scale instances. In the rest of this section, our discussion departs from~\cite{ccopf-sirev}, highlighting a number of enhancements over the algorithm they propose.}

The cutting-plane approach in \cite{ccopf-sirev} iteratively solves a sequence of relaxations of~\eqref{eq:robustprob}. At each iteration, we must verify if~\eqref{eq:robustprob} is satisfied. \textcolor{black}{In the case of partitioned uncertainty,}~\eqref{eq:robustprob} holds iff
\begin{equation}\label{eq:max}
\left[\max_{\mu \in U_\mu} x^T\mu\right] + \Phi^{-1}(1-\epsilon)\sqrt{\max_{\Sigma \in U_\Sigma} x^T\Sigma x} \leq b.
\end{equation}

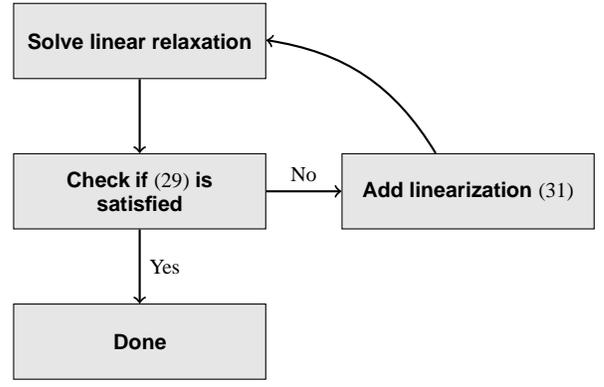
\begin{figure}[t]

\centering
\begin{tikzpicture}
[node distance = 1cm, auto,font=\footnotesize,
every node/.style={node distance=2cm},
force/.style={rectangle, draw, fill=black!10, inner sep=5pt, text width=3cm, text badly centered, minimum height=1cm, font=\bfseries\footnotesize\sffamily}] 

\node [force] (check) {Check if~\eqref{eq:robustprob} is satisfied};
\node [force, above of=check] (solve) {Solve linear relaxation};
\node [force, right=1cm of check] (add) {Add linearization~\eqref{eq:robustlinearization} };
\node [force, below of=check] (done) {Done};


\path[->,thick] 
(add) edge [bend right=25] (solve);

\draw [->,thick] (solve) -- (check);
\draw [->,thick] (check) -- (done) node [midway] {Yes};
\draw [->,thick] (check) -- (add) node [midway] {No};

\end{tikzpicture} 
\caption{\textcolor{black}{An illustration of the iterative cutting-plane approach used to enforce the convex, nonlinear distributionally robust chance constraints.}}
\label{fig:cp}
\end{figure}

\noindent For fixed $x^*$, both of the inner maximization problems in~\eqref{eq:max} have a linear objective, and so detecting if~\eqref{eq:robustprob} is satisfied for a given $x^*$ can be computed by optimizing a linear function over the sets $U_\mu$ and $U_\Sigma$. 
If the solution to \eqref{eq:max} shows that \eqref{eq:robustprob} is satisfied, then the algorithm terminates. Otherwise, \eqref{eq:max} is used to find the corresponding $(\mu^*,\Sigma^*)$ that violates \eqref{eq:robustprob}, and we add a \textit{linearization} of the corresponding constraint~\eqref{eq:reform} \textcolor{black}{of the form}
\begin{gather}
\color{black}
x^T\mu^* + \Phi^{-1}(1-\epsilon)\sqrt{(x^*)^T\Sigma^*x^*} + \nonumber \\
\color{black}
\left(\Phi^{-1}(1-\epsilon)/\sqrt{(x^*)^T\Sigma^*x^*}\right)(x^*)^T\Sigma^*(x-x^*) \le b
\label{eq:robustlinearization}
\end{gather}
to the relaxation which cuts off the current solution. Therefore, at any iteration, the relaxation we solve is a linear program, similar to the linearization scheme for CC OPF. This process repeats until~\eqref{eq:robustprob} is satisfied within numerical tolerances.
\textcolor{black}{Our approach differs in two notable aspects from the algorithm proposed by \cite{ccopf-sirev}:
    \begin{enumerate}
        \item They propose to treat the term $\max_{\mu \in U_\mu} x^T\mu$ by standard reformulation techniques based on strong duality to generate an equivalent deterministic formulation in an extended set of variables. Instead, based on the empirical observation that only a few extreme cases are important, we also apply the cutting-plane technique to this term. In other cases, it may be advantegous to reformulate~\cite{cutreform}.
        \item They propose to introduce a slack variable $t$ for $\sqrt{\max_{\Sigma \in U_\Sigma} x^T\Sigma x}$ and add a linearization whenever the constraint $t \ge \sqrt{\max_{\Sigma \in U_\Sigma} x^T\Sigma x}$ is violated. However, it is possible for this constraint to be violated when the original~\eqref{eq:robustprob} is not. Our implementation avoids unnecessary iterations and only adds linearizations when~\eqref{eq:robustprob} is violated. 
\end{enumerate}
}
Although this algorithm does not have polynomial convergence guarantees in general, it is an immensely powerful approach that also underlies standard techniques such as Benders decomposition. \textcolor{black}{See Figure~\ref{fig:cp} for an illustration of the algorithm. We refer readers to \cite{ccopf-sirev} for further discussion.}

\color{black}
We also note an enhancement for this algorithm based on the special structure of the power output and ramping constraints~\eqref{eq:CC-OPF_pmax}-\eqref{eq:CC-OPF_RD}.
For example, consider the constraint
\begin{equation}\label{eq:linearchance}
 P(p_i - \boldsymbol\Omega \alpha_i \le p_i^{max}) \ge 1-\epsilon_i.
\end{equation}
Recall that $\boldsymbol\Omega = \sum_{b \in \mathcal{W}} \boldsymbol\omega_b$ is the total deviation from the forecast, hence $\boldsymbol\Omega$ is a scalar Gaussian random variable with mean $\sum_{b \in \mathcal{W}} \mu_b$ and variance $e^T\Sigma e$ where $e$ is the vector of all ones.
Hence~\eqref{eq:linearchance} is equivalent to
\begin{equation}
p_i-\left(\sum_{b \in \mathcal{W}} \mu_b\right)\alpha_i + \Phi^{-1}(1-\epsilon)\sqrt{\alpha_i^2 e^T\Sigma e } \le p_i^{max}.
\end{equation}

In the CC-OPF case where we assume zero means and $\alpha_i \ge 0$, this constraint simplifies further to
\begin{equation}
p_i + \Phi^{-1}(1-\epsilon)\alpha_i\sqrt{ e^T\Sigma e } \le p_i^{max},
\end{equation}
which is a single linear constraint; no SOC constraints are needed, an important property which was not noted by \cite{ccopf-sirev}. Furthermore, in the distributionally robust case, the worst-case realizations may be computed independently of $\alpha_i$, i.e.,
\begin{equation}\label{eq:robustequiv}
p_i - \alpha_i \left[\min_{\mu \in U_\mu} \left(\sum_{b \in \mathcal{W}} \mu_b\right)\right] +\Phi^{-1}(1-\epsilon)\alpha_i\sqrt{\max_{\Sigma \in U_\Sigma} e^T\Sigma e } \le p_i^{max}.
\end{equation}
The constraint~\eqref{eq:robustequiv} is again a single linear constraint. This observation is useful for constraints~\eqref{eq:CC-OPF_pmax}-\eqref{eq:CC-OPF_RD} with the structure where a single decision variable multiplies a scalar Gaussian random variable, since we may avoid applying the cutting-plane algorithm to these constraints. This special structure does not hold for the line capacity constraints~\eqref{eq:CC-OPF_fplus}-\eqref{eq:CC-OPF_fminus}.
\color{black}

\textcolor{black}{Finally, we note that for a small number of special cases, including when $U$ cannot be partitioned, Ben-Tal et al~\cite{ben2009robust} derive explicit semidefinite programming representations of~\eqref{eq:robustprob}. We leave an exploration of this case for future work.  }

\subsection{Formulation for RCC OPF}\label{sec:RCC-OPF_formulation}

We adapt the discussion in \textcolor{black}{Section}~\ref{sec:robustmodel} to formulate an RCC OPF. First, we adopt the assumptions in \cite{ccopf-sirev} that the fluctuations at the different wind sites are independent within an OPF time step, i.e. $\Sigma = \operatorname{diag}(\sigma^2)$. \textcolor{black}{For brevity, we denote the multivariate Gaussian distribution with diagonal covariance matrix as $N(\mu,\sigma^2)$ where $\sigma^2$ is a vector.} The RCC OPF models uncertainty in the $\boldsymbol\omega$ distribution parameters at each bus $b$ as intervals $[-\bar \mu_b,\bar\mu_b]$ and $[\sigma_b^2 - \bar\sigma_b^2,\sigma_b^2 + \bar\sigma_b^2]$ for $\mu$ and $\sigma$, respectively. To represent the aggregate uncertainty, we \textcolor{black}{follow~\cite{ccopf-sirev} and} construct polyhedral uncertainty sets
\begin{equation}\label{eq:U_mu}
U_\mu = \left\{ \mu \in \mathbb{R}^{|\mathcal{W}|} : |\mu_b| \leq \bar\mu_b, \sum_{b \in \mathcal{W}} \frac{|\mu_b|}{\bar\mu_b} \leq \Gamma_\mu|\mathcal{W}| \right\}
\end{equation}
and
\begin{equation}\label{eq:U_sigma_sq}
\begin{split}
\color{black} U_{\sigma^2} = \Bigg\{ s \in \mathbb{R}^{|\mathcal{W}|} : \exists t: s = \sigma^2 + t,  \\
\color{black} |t_b| \leq \bar\sigma_b^2, \sum_{b \in \mathcal{W}} \frac{|t_b|}{\bar\sigma_b^2} \leq \Gamma_\sigma|\mathcal{W}| \Bigg\}
\end{split}
\end{equation}

\noindent for the mean and variance, similar to those proposed by Bertsimas and Sim~\cite{Bertsimas2004}.

\textcolor{black}{In our RCC OPF formulation, we replace constraints~\eqref{eq:CC-OPF_pmax}-\eqref{eq:CC-OPF_fminus} with their distributionally robust counterparts.}
\begin{gather}
\color{black}
 P_{\boldsymbol\omega \sim N(\mu,\sigma^2)}(p_i - \boldsymbol\Omega \alpha_i \le p_i^{max}) \ge 1- \epsilon_i, \nonumber\\
 \color{black}
 \forall i \in \mathcal{G}, \forall \mu \in U_\mu, \sigma^2 \in U_{\sigma^2}\label{eq:CC-OPF_pmax_robust} \\
 \color{black}
 P_{\boldsymbol\omega \sim N(\mu,\sigma^2)}(p_i - \boldsymbol\Omega \alpha_i \ge p_i^{min}) \ge 1- \epsilon_i, \nonumber\\
 \color{black}
 \forall i \in \mathcal{G}, \mu \in U_\mu, \sigma^2 \in U_{\sigma^2}\label{eq:CC-OPF_pmin_robust} \\
 \color{black}
 P_{\boldsymbol\omega \sim N(\mu,\sigma^2)}(-\boldsymbol\Omega\alpha_i \le RU_i) \ge 1-\epsilon_i, \nonumber\\
 \color{black}
 \forall i \in \mathcal{G}, \mu \in U_\mu, \sigma^2 \in U_{\sigma^2}\label{eq:CC-OPF_RU_robust} \\
 \color{black}
 P_{\boldsymbol\omega \sim N(\mu,\sigma^2)}(\boldsymbol\Omega\alpha_i  \le RD_i) \ge 1- \epsilon_i, \nonumber\\
 \color{black}
 \forall i \in \mathcal{G}, \mu \in U_\mu, \sigma^2 \in U_{\sigma^2}\label{eq:CC-OPF_RD_robust} \\
 \color{black}
P_{\boldsymbol\omega \sim N(\mu,\sigma^2)}\bigg(f_{mn} + \beta_{mn}\boldsymbol\Omega(\delta_n - \delta_m) + \nonumber\\
\color{black}
\beta_{mn}\boldsymbol\omega^T(\pi_m-\pi_n) \le f_{mn}^{max} \bigg)  \ge 1- \epsilon_{mn},\nonumber\\
\color{black}
\quad \forall \{m,n\} \in \mathcal{L},  \mu \in U_\mu, \sigma^2 \in U_{\sigma^2}\label{eq:CC-OPF_fplus_robust} \\
\color{black}
P_{\boldsymbol\omega \sim N(\mu,\sigma^2)}\bigg( f_{mn} + \beta_{mn}\boldsymbol\Omega(\delta_n - \delta_m) +\nonumber\\
\color{black}
\beta_{mn}\boldsymbol\omega^T(\pi_m-\pi_n) \ge -f_{mn}^{max} \bigg) \ge 1- \epsilon_{mn}, \nonumber\\
\color{black}
\forall \{m,n\} \in \mathcal{L},  \mu \in U_\mu, \sigma^2 \in U_{\sigma^2}\label{eq:CC-OPF_fminus_robust}
\end{gather}

\textcolor{black}{For the objective~\eqref{eq:ccobjective}, we simply take the nominal value for $\operatorname{var}(\boldsymbol\Omega)$, leaving this term unchanged. }
The parameters $\Gamma_\mu$ and $\Gamma_\sigma$ are the \textit{uncertainty budgets} used to adjust the level of conservatism of the resulting RCC OPF algorithm. The $\Gamma$'s may be interpreted as a bound on the proportion of wind farms that may take their \textit{worst case distribution}. The least conservative limit $\Gamma = 0$ recovers the standard CC OPF model. The most conservative limit $\Gamma = 1$ ensures feasibility when each wind farm can take on its worst-case production distribution. \textcolor{black}{Uncertainty budgets have a well-studied probabilistic interpretation for the case of robust \textit{linear} constraints~\cite{Bertsimas2004}.  In this work, however, we enforce robust feasibility of a second-order cone constraint with respect to perturbations in the coefficients. One \textit{could} interpret our formulation as presuming some prior distribution on the parameters of the distribution of $\boldsymbol\omega$ and then requiring feasibility of the chance constraint under high probability in the space of distributions. However, few probabilistic guarantees are known for this case; see~\cite[Ch. 10]{ben2009robust}.}

\textcolor{black}{
Our motivation for choosing this formulation, although it lacks a rigorous probabilistic interpretation in terms of prior distributions, is that it adds very little computational expense to the CC OPF problem. Indeed, for the uncertainty sets in \eqref{eq:U_mu} and \eqref{eq:U_sigma_sq}, the optimization \eqref{eq:max} that required to check feasibility of each robust chance constraint~\eqref{eq:robustprob} reduces to a cheap sorting operation~\cite{Bertsimas2004,cutreform}. At the same time, as we intend to demonstrate in the case study, the RCC OPF formulation delivers measurable benefits in practice.}


Note that, similar to the CC OPF, the constraints \eqref{eq:CC-OPF_pmax_robust}-\eqref{eq:CC-OPF_fminus_robust} are treated separately, i.e. when they are ``robustified'' with respect to the distribution parameters, the worst-case distribution may be different for each constraint. In the remainder of this manuscript, we set $\Gamma = \Gamma_\mu = \Gamma_\sigma$ and investigate the results as a function of $\Gamma$.

\textcolor{black}{We note a connection to the robust OPF models discussed in in \cite{guan_2013,lorca_2014,bertsimas_2013,jabr_2014, rjabr_2013b}. The reformulation from ~\eqref{eq:reform} to ~\eqref{eq:prob} under the assumption of Gaussian uncertainty has a 1-1 equivalence with robust linear constraints with ellipsoidal uncertainty sets~\cite{bertsimas_2004}. Hence, our CC OPF may be interpreted as a particular form of robust OPF, and under that interpretation, the effect of the RCC OPF formulation is to introduce a level of uncertainty sets on the nominal values of the robust problem.}
 

\section{Case Study}\label{sec:case_study}

We investigate the performance and benefits of an RCC OPF approach relative to CC OPF and deterministic OPF approaches by implementing all three in the setting of vertically-integrated grid operations, specifically, on a modification of the Bonneville Power Administration (BPA) system\footnote{\textcolor{black}{We refer interested readers to \cite{ccopf-sirev} for illustrative case studies of the CC OPF model on the IEEE Realibility Test System.}}.  This setting allows us to evaluate the effect of using probabilistic methods in an OPF. We use $\Gamma$ as a parameter to study the impact of robust conservatism on the cost of generation and the statistics of Area Control Error (ACE), generator ramping, and power line loading.  

\subsection{Test System and Data}\label{sec:test_system}
We use a modification of the BPA system\footnote{\textcolor{black}{For visual representations of the BPA system and some of its partitions, interested readers are referred to \cite{bpa_map}.}} with 2866 transmission lines and 2209 buses including 676 load buses, 176 controllable generators, and 24 wind farms. The total installed capacity of controllable generators is 40.6 GW, composed of 133 hydro generators (28 GW), 41 gas-fired generators (9.6 GW), 1 coal-fired generator (1.2 GW), and 1 nuclear generator (1.8 GW). The total installed capacity of 24 wind farms is 4.6 GW. The technical characteristics of the controllable generators and the network configuration are adapted from PowerWorld\cite{powerworld}. Tie-line power flows to neighboring interconnections are modeled as loads at the ends of the tie lines. 

The case study spans the Winter season from December 2012 to March 2013 so that the factors affecting wind, load, and hydro generation are relatively stationary.  Hour-ahead, hour-resolution load forecasts and their actual 5-minute realizations are taken from \cite{bpa}. As provided, these data are aggregated across the entire BPA system. We disaggregate the load forecasts and realizations among 676 load buses based on their population density. For each wind farm location, archived NOAA forecast data \cite{noaa} are used to provide hour-ahead, hour-resolution wind speed forecasts. Actual, 5-minute resolution wind speed realizations at each wind farm location are taken from BPA historical data \cite{bpa}. Both the forecasts and realizations are converted to wind power using an equivalent wind turbine for each wind farm location \cite{hayes_2011}. We use the methodology from \cite{letter} to compute the intervals  $[-\bar \mu_b,\bar\mu_b]$ and $[\sigma_b^2 - \bar\sigma_b^2,\sigma_b^2 + \bar\sigma_b^2]$ for the uncertainty sets in Eq. \eqref{eq:U_mu} and \eqref{eq:U_sigma_sq}. \textcolor{black}{The methodology in \cite{letter} uses a data-driven statistical analysis, which reveals that hourly average wind speed is proportional to the intra-hour wind speed variability. Next, this relationship is used to fit historical wind speed forecast errors to a generalized normal distribution. This fitting yields the best-fit ranges on the hourly-average wind speed and its standard deviation. After that, these ranges are converted from wind speed to wind power units by using a typical wind turbine power curve.}

In this study, we seek to realistically emulate the short-term operational planning of the BPA system operator; however, not all required data are available. To fill in the missing data, we make the following operational assumptions:
\paragraph{Nuclear Generators} The single nuclear unit in the BPA system is assumed to be a ``must run'' unit, and its hourly power output $p_i$ is set to the 95\% of its nameplate capacity. Its participation factor $\alpha_i$ in real-time balancing is set to zero \cite{jaleeli_1992}. 
\paragraph{Gas and Coal Generators} The power outputs $p_i$ and participation factors $\alpha_i$ of all gas-fired and one coal-fired generators are decision variables.
\paragraph{Hydro Generators} Dispatch decisions for hydro generators often depend more on water flow considerations rather than on power system operations \cite{piekutowski_1994}. Instead of being co-optimized with thermal and nuclear generators, hydro dispatch levels $p_i$ are fixed in all OPF formulations considered here  as exogenous parameters \cite{studarus_2013}. System aggregated, hourly-resolution hydro generation is taken from historical BPA data \cite{bpa} and disaggregated to individual hydro generators based on their installed capacity.  The assignment of participation factors is also affected by water flow conditions beyond the scope of this work. Lacking operational data, we set the participation factors $\alpha_i$ of the hydro generators to a common value, which itself is a decision variable. 

We refer interested readers to \cite{release}  for the input data and the code used in this case study. 

\subsection{Evaluation Procedure} \label{sec:evaluation_procedure}
The evaluation procedure includes two steps, which emulate the hour-ahead scheduling and real-time dispatch, respectively, and are organized as follows:
\subsubsection{Step 1} The RCC OPF, CC OPF (i.e. RCC OPF with $\Gamma$=0), and a deterministic OPF are solved for hour $t$ using hour-ahead, hour-resolution wind power forecasts $\omega^f_{t,b}$ and load forecasts $d^f_{t,b}$ at each bus $b$ and the generator commitment decisions $u_{t,i}$ described in Section \ref{sec:test_system}. The result is an optimal hourly dispatch $p^*_{t,i}$ and hourly participation factors $\alpha^*_{t,i}$. 
\subsubsection{Step 2}  Next, these optimal decisions $p^*_{t,i}$ and $\alpha^*_{t,i}$ are evaluated in a quasi-static power flow simulation of the system behavior using the actual, 5-minute realizations of wind power $\omega_{t,b}(\tau)$ and demand $d_{t,b}(\tau)$ where $\tau$ refers to the twelve 5-minute intervals of hour $t$. For every $\tau$, we compute the actual power output of each controllable generator $p_{t,i}(\tau)$ as $p_{t,i}(\tau) = \max[p^{min}_i,\min[\hat p_{t,i}(\tau),p^{max}_{i}]]$ where
\begin{equation}
\hat p_{t,i}(\tau)=p^*_{t,i} - \alpha^*_{t,i} \sum_{j \in \mathcal{W}} \omega_{t,b}(\tau).
\end{equation}
Using $p_{t,i}(\tau)$, a DC power flow calculation yields the $f_{t,mn}(\tau)$. The actual power output $p_{t,i}(\tau) $ is then used to calculate the actual hourly operating cost $C_t$ by summing the cost of $p_{t,i}(\tau)$ over all 5-minute intervals $\tau$. This emulation process reflects the vertically-integrated utility setting, i.e. costs are allocated according to energy delivered to the system with no markup cost for providing regulation.

Steps 1 and 2 are repeated for every operating hour in the period from December 2012 to March 2013. The results for each $\tau$ are analyzed for power system area control error (ACE) statistics, generator ramping statistics, and power line flow statistics as described below.

\subsection{\textcolor{black}{Cost Performance}}

The RCC OPF includes three user-defined parameters related to its probabilistic nature.  The first of these is $\Gamma$ which determines the budget of uncertainty in wind forecast probability distributions defined by the uncertainty sets in \eqref{eq:U_mu} and \eqref{eq:U_sigma_sq}. Wind conditions change frequently, and we expect that this uncertainty parameter will be determined by short-term policy decisions of the system operators.  In contrast, the parameters $\epsilon_i$ and $\epsilon_{mn}$ limit the probability that equipment constraints are violated, i.e. generation output or ramping limits described in \eqref{eq:CC-OPF_pmax}-\eqref{eq:CC-OPF_RD} and power flow limits on lines described in \eqref{eq:CC-OPF_fplus}-\eqref{eq:CC-OPF_fminus}. These parameters are directly related to the impact on power system assets, and we expect these  are determined by long-term, e.g. seasonal policies. 

In our case study, these parameters are set in a sequential process with the results shown in Tables \ref{table:cost_gamma} and \ref{table:cost_epsilon}. Fixing the parameters $\epsilon_i$ and $\epsilon_{mn}$ at $1/6$ and 0.0025, respectively, Steps 1 and 2 from \textcolor{black}{Section}~\ref{sec:evaluation_procedure} are executed to determine the actual cost of generation $C^a$ over the entire study period for $\Gamma$ between 0.0\footnote{Recall that the CC OPF is the RCC OPF with $\Gamma$=0.} and 1.0 and for a deterministic OPF. The participation factors are not decision variables in the deterministic OPF. Instead, a fixed participation factor $\alpha$ = 0.05 is used for all thermal and hydro generators. Table \ref{table:cost_gamma} shows the lowest $C^a$ is found at $\Gamma^*$=0.6. 

Fixing $\Gamma=\Gamma^*=0.6$, Table~\ref{table:cost_epsilon} displays the sensitivity of the actual cost $C^a$ of the RCC OPF solution to the parameters $\epsilon_i$ and $\epsilon_{mn}$.  We begin the discussion with the results in the row for $\epsilon_{mn}$=0.01.  For this larger $\epsilon_{mn}$, potential network congestion plays less of a role, i.e. the selection of the generators' $p_{t,i}^*$ and $\alpha^*_{t,i}$ is less dependent on their location in the network. Instead, their selection is more sensitive to generator costs and constraints. Therefore, as $\epsilon_i$ decreases (moving left to right) in the $\epsilon_{mn}$=0.01 row, the main impact is to spread the $p_{t,i}^*$ and $\alpha^*_{t,i}$ more uniformly across the fleet to reduce the ramping duty of any particular generator.  As $\epsilon_i\rightarrow 1/48$, more duty is placed on higher cost generators driving up the operating cost via the $\operatorname{var}(\boldsymbol\Omega)$ term in \eqref{eq:CC-OPF_objective}. However for our case study, the high percentage of very flexible hydro generators suppresses this cost increase---a result that is not expected to carry over to other power systems with different controllable generation fleets. 

Next, we consider the small $\epsilon_{mn}$ limit shown in the $\epsilon_{mn}$=0.0001 row in Table~\ref{table:cost_epsilon}. Here, avoiding potential network congestion plays a larger role in the selection of the $p_{t,i}^*$ and $\alpha^*_{t,i}$ with the results becoming less sensitive to the generator cost and risk parameter $\epsilon_i$. This is reflected in the elevated and relatively flat cost even as $\epsilon_i\rightarrow 1/48$. In between the two extremes of $\epsilon_{mn}$, there is a relatively strong trade off in $C^a$ between $\epsilon_{mn}$ (network risk) and $\epsilon_i$ (generator risk). 

The results in Tables \ref{table:cost_gamma} and \ref{table:cost_epsilon} suggest that, for the test system used in our case study, the RCC OPF model achieves the best cost performance with $\Gamma^* = 0.6$, $\epsilon_{mn}^* = 0.0025$, $\epsilon_i^* = 1/6$. In the remainder of this manuscript, we assume that the long-term policy parameters are fixed at $\epsilon_{mn}^* = 0.0025$ and $\epsilon_i^* = 1/6$,  and present the technical analysis for variable $\Gamma$.

\begin{table}[htbp!]
\caption{Cost performance of the OPF models in the period from December 2012 to March 2013}
\begin{center}
\label{table:cost_gamma}
\begin{tabular}{>{\centering\arraybackslash}m{0.3cm} | >{\centering\arraybackslash}m{0.6cm}| >{\centering\arraybackslash}m{0.8cm} >{\centering\arraybackslash}m{0.8cm} >{\centering\arraybackslash}m{0.9cm} >{\centering\arraybackslash}m{0.4cm} >{\centering\arraybackslash}m{0.5cm}| >{\centering\arraybackslash}m{0.4cm} }
\hline\hline
  &  CC OPF   & \multicolumn{5}{c|}{RCC OPF} & OPF\\
\hline
  $\Gamma$  &  0 & 0.2 & 0.4 &  $\mathbf{0.6}$ & 0.8 & 1.0 & --\\
\hline\hline
   $C^{a}$, M\$ & 112.2 & 111.6 & 110.9&  $\mathbf{108.7}$ & 112.7 & 114.9& 115.6\\
\hline
    $\Delta$, M\$ & -- & -0.589 & -1.296 &  -$\mathbf{3.471}$ & 0.459 & 2.716& 3.401\\
\hline
    $\Delta$, \% & -- & -0.524 & -1.155 &  -$\mathbf{3.093}$ & 0.409 & 2.420& 3.032\\
 \hline\hline
\multicolumn{8}{p{0.9\linewidth}}{The actual realized generation cost $C^a$ for the period from December 2012 to March 2013 is computed in Step 2 using the dispatches $p_{t,i}^*$ and participation factors $\alpha_{t,i}^*$ computed in Step 1 of \textcolor{black}{Section}~\ref{sec:evaluation_procedure}. The least cost solution is found for $\Gamma$=0.6 and is marked in bold. Also displayed are the changes in cost and fractional changes in cost relative to the $\Gamma$=0 case.    In all of these cases, $\epsilon_{mn}$=0.0025 and $\epsilon_i$= and $1/6$. }
\end{tabular}
\end{center}
\end{table}

\begin{table}[htbp!]
\caption{Sensitivity of the RCC OPF cost (\%*) to $\epsilon_i$ and $\epsilon_{mn}$}
\begin{center}
\begin{tabular}{ c| c c c c }
\hline\hline
  & $\epsilon_i = \frac{1}{6}$   & $\epsilon_i = \frac{1}{12}$ & $\epsilon_i = \frac{1}{24}$ &  $\epsilon_i = \frac{1}{48}$ \\
\hline\hline
$\epsilon_{mn} = 0.01$  & 1.379 & 1.471 & 1.563 &  1.563 \\
\hline
$\epsilon_{mn} = 0.005$  & 0.827& 1.195 & 1.379 &  1.471 \\
\hline
$\epsilon_{mn} = 0.0025$  & $\mathbf{0}$ & 1.379 & 1.379 &  1.471 \\
\hline
$\epsilon_{mn} = 0.001$  & 0.460 & 1.011 & 1.471 &  1.471 \\
\hline
$\epsilon_{mn} = 0.0001$ & 2.849 & 2.941 & 2.941 &  2.941 \\
 \hline\hline

\multicolumn{5}{p{0.8\linewidth}}{Percentage changes in actual generation cost $C^a$ relative to the $\epsilon_i=1/6$, $\epsilon_{mn}=0.0025$ case.  All cases use $\Gamma=\Gamma^*=0.6$.}
\end{tabular}

\end{center}
\label{table:cost_epsilon}
\end{table}

\subsection{\textcolor{black}{Technical Analysis}}
\subsubsection{\textcolor{black}{ACE Performance}}
The Area Control Error (ACE) is computed for each 5-minute interval $\tau$ in hour $t$ as:
\begin{equation}
ACE_t(\tau) = \sum_{b \in \mathcal{B}} (d_{t,b} (\tau) - w_{t,b} (\tau)) - \sum_{i \in \mathcal{G} } p_{t,i}(\tau) . 
\end{equation}
\begin{figure}[b!]
  \centering
    \includegraphics[width=\linewidth]{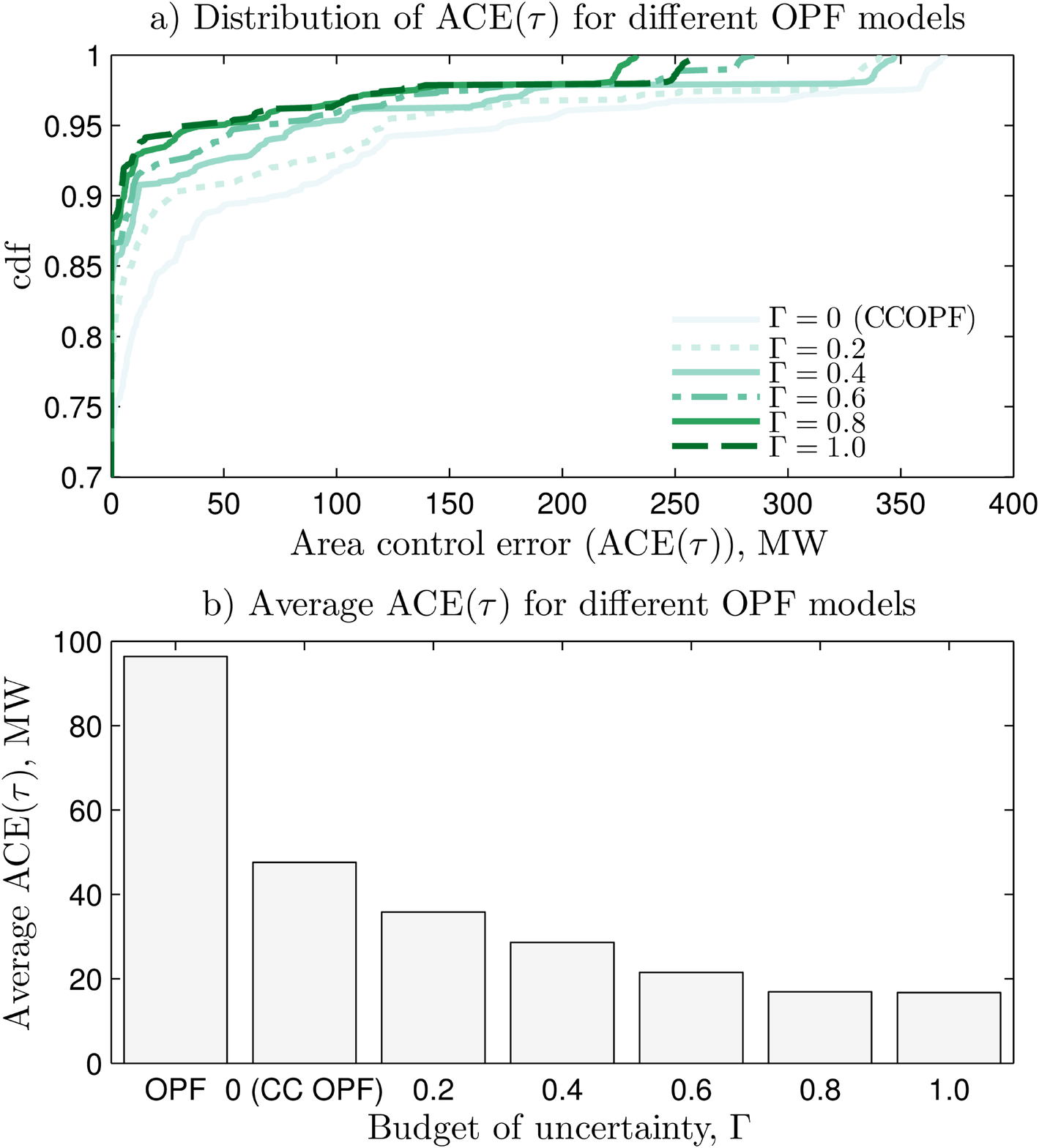}
\caption{a) Cumulative distribution function (CDF) of $ACE_t(\tau)$ for different values of $\Gamma$=0.0 to 1.0. The CDF for the deterministic OPF is not shown for clarity of the Figure. b) Average of $ACE_t(\tau)$ for the same values of $\Gamma$ as in a). Under the conditions of the study, the generators always have sufficient downward flexibility to avoid $ACE_t(\tau) <0$, i.e. overgeneration is not observed. On the other hand, undergeneration ($ACE_t(\tau) >0$) occurs for all models. All results are computed with $\epsilon_{mn}^* = 0.0025$ and $\epsilon_i^* = 1/6$. }
\label{fig:ace}
\end{figure}
Figure \ref{fig:ace}a displays the cumulative distribution function (CDF) of all the $ACE_t(\tau)$ in the study period, and \textcolor{black}{Figure}~\ref{fig:ace}b displays the average of $ACE_t(\tau)$. Starting from the most conservative $\Gamma$ = 1.0, the average of $ACE_t(\tau)$ displays a slow but monotonic increase showing that $\Gamma$ effectively controls the system's technical performance. At $\Gamma <$ 0.6, the average $ACE_t(\tau)$ increases more rapidly to the CC OPF at $\Gamma$=0 and displays a significant jump from the CC OPF to the deterministic OPF. By accounting for fluctuations in wind, the CC OPF outperforms the deterministic OPF in controlling the ACE, and by progressively accounting for uncertainty in the parameters of the distribution describing the fluctuations, the RCC OPF outperforms the CC OPF. 

In \textcolor{black}{Figure}~\ref{fig:ace}a), the difference between the CDFs for $\Gamma$ = 1.0, 0.8, and 0.6 is not very significant. As our measure of conservatism is relaxed further (i.e. $\Gamma <$0.6), the CDFs show a general increase in the frequency of ACE events of all sizes and the emergence of a longer tail of large $ACE_t(\tau)$ values. We also note that $\Gamma \sim$ 0.6 is the value where the ACE statistics first begin to significantly deteriorate and where RCC OPF cost takes on its minimum value. Above $\Gamma =$0.6, little additional ACE control performance is gained for the additional cost. This analysis suggests that, in the setting of \textcolor{black}{vertically-integrated} grid operations, the RCC OPF with an appropriately chosen $\Gamma$ will result in better compliance with the control performance standards (CPS) \cite{cps} at a lower operating cost. 
\subsubsection{\textcolor{black}{Ramping Performance}}
\begin{figure}[htbp!]
  \centering
    \includegraphics[width=\linewidth]{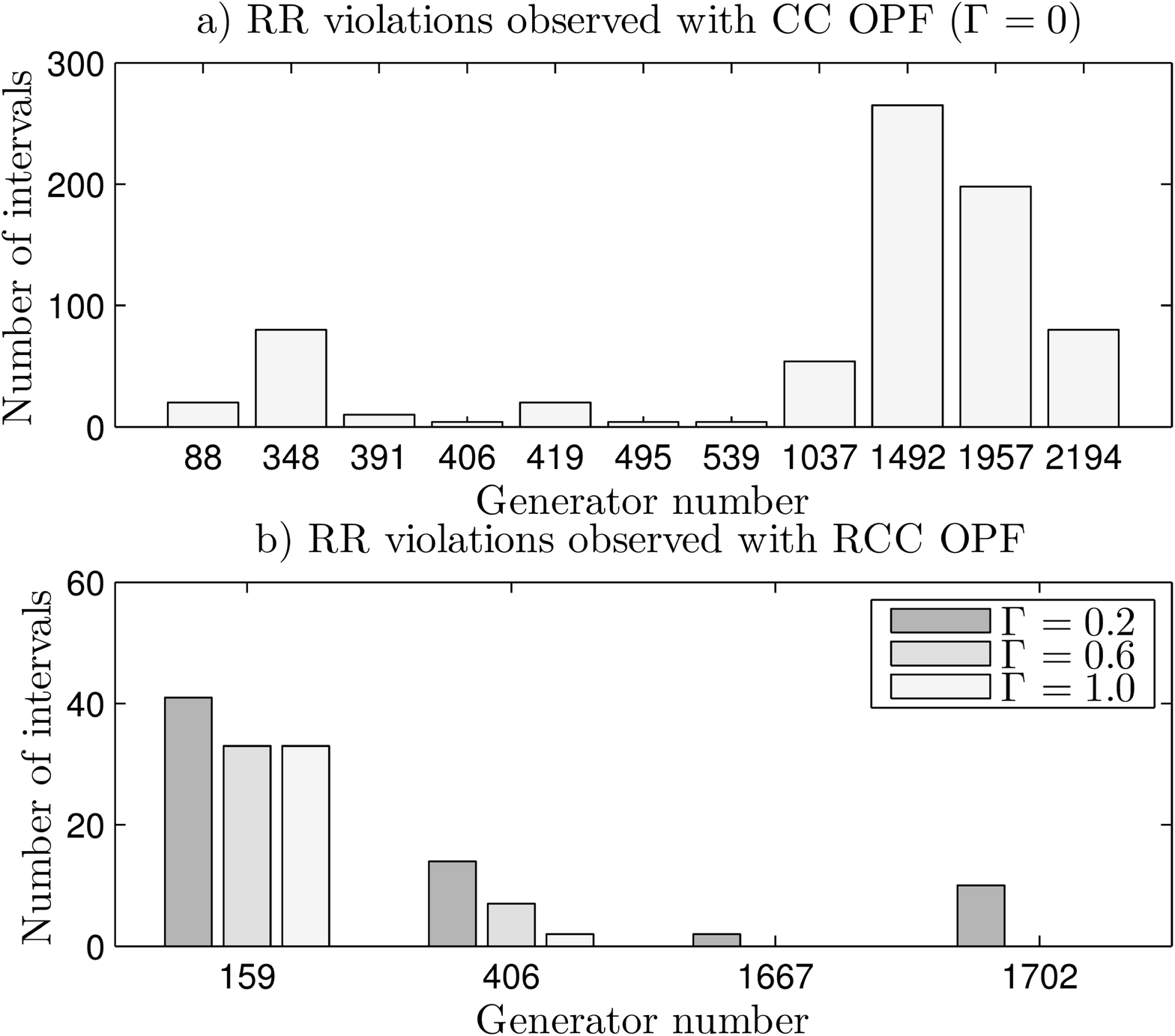}
\caption{a) Histogram of the number of generator ramp rate (RR) violations per generator over the entire study period for the CC OPF (i.e. the RCC OPF with $\Gamma$=0). b) Same as a) but for the RCC OPF with $\Gamma$ = 0.2, 0.6, and 1.0. All results are computed with $\epsilon_{mn}^* = 0.0025$ and $\epsilon_i^* = 1/6$.}
\label{fig:ramp}
\end{figure}

\indent The RCC OPF also reduces generator ramp rate (RR) violations as compared to the CC OPF and deterministic OPF, potentially avoiding generator wear-and-tear effects \cite{troy_2012}. After the $\alpha_{t,i}^*$ are chosen, the generator ramp rates are simply $\alpha_{t,i}^* \cdot \sum_{b \in \mathcal{W}} \omega_{t,b}(\tau)$. Figure~\ref{fig:ramp}a) and b) display the number of RR violations for individual generators for the CC OPF and for the RCC OPF for different $\Gamma$, respectively.  

The impact of the RCC OPF on RR violations is twofold. First, as the robustness of the RCC OPF increases, the number of generators affected by RR violations is reduced from 11 (for the CC OPF, i.e. $\Gamma$=0), to 4 with $\Gamma$=0.2, and to 2 with $\Gamma$=1.0. Second, the number RR violations per generator is also greatly reduced as $\Gamma$ increases. It is noteworthy that both effects can be observed even for a relatively small level of robustness, e.g. $\Gamma = 0.2$. \textcolor{black}{As shown in Figure~\ref{fig:ramp}b), the greatest number of 5-minute intervals when the robust chance constraints on ramping are violated, is observed with $\Gamma = 0.2$ for generator 159. This translates into the fraction $\approx$0.004 of the total number of intervals considered in this case study, which is less than $\epsilon^*_i=1/6$ enforced in the chance constraints.} Combined with the results from Table \ref{table:cost_gamma}, this analysis suggests that the RCC OPF with an appropriately chosen $\Gamma$ achieves lower operating cost and avoids indirect costs related to wear-and-tear effects on controllable generators. In contrast, RR violations for the deterministic OPF (not shown in \textcolor{black}{Figure}~\ref{fig:ramp}) are observed on 39 generators---24 of which experience RR violations in more than one 5-minute interval. 
\subsubsection{\textcolor{black}{Transmission Overload Performance}}
\begin{figure}[t!]
  \centering
    \includegraphics[width=\linewidth]{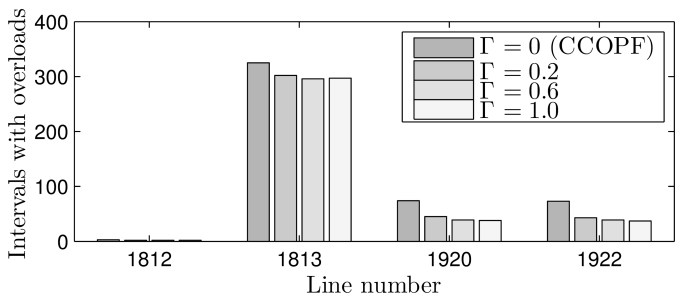}
\caption{Comparison of overloads that are observed with the CC and RCC OPF models during more than one 5-minute interval.  All results are computed with $\epsilon_{mn}^* = 0.0025$ and $\epsilon_i^* = 1/6$.  
}
\label{fig:overload}
\end{figure}

From the $p_{t,i}(\tau)$, $\omega_{t,b}(\tau)$, and $d_{t,b}(\tau)$, a power flow solution yields $f_{t,mn}(\tau)$ from which power line overloads are computed. Figure~\ref{fig:overload} displays a histogram of the number of overloads per power line for the four most frequently overloaded lines in the RCC OPF for $\Gamma$ = 0.0, 0.2, 0.6, and 1.0. Several other lines are overloaded during the study period, but these overloads only occur during one 5-minute period.  Interestingly, $\Gamma$ does not have a significant impact on the frequency of overloads for the lines in \textcolor{black}{Figure}~\ref{fig:overload}.

\textcolor{black}{Among the most overloaded lines in Figure \ref{fig:overload}, the greatest number of violations of robust chance constraints on power flow limits is observed on line 1813 with $\Gamma = 0$. These violations are observed on the fraction $\approx$0.013 of empirical realizations, which is larger than $\epsilon_{ij}^{*}=0.0025$ enforced in the chance constraints\footnote{\textcolor{black}{The exact cause of violations on line 1813 cannot easily be explained by particular attributes of the model or the BPA test system, but the issue can be resolved via out-of-optimization corrections \cite{studarus_2013,hedman_2015} used in practice.}}. For other lines in Figure \ref{fig:overload}, the empirical number of violations are less than $\epsilon_{ij}^{*}=0.0025$. }

\subsection{\textcolor{black}{Computational Performance}}
\begin{figure}[t!]
  \centering
    \includegraphics[width=\linewidth]{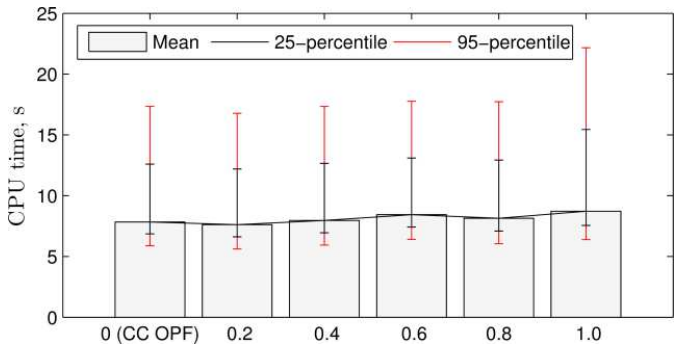}
\caption{\textcolor{black}{Computational performance of an RCC OPF instance.}} \label{fig:cpu_time}
\end{figure}

The computations were carried out with CPLEX 12.6 \cite{cplex} as an LP solver on a Intel Xenon 2.55 GHz processor with at least 32 GB RAM on the Hyak supercomputer system at the University of Washington \cite{hyak}. All modeling was done using JuMPChance~\cite{jumpchance}, a freely available extension for the JuMP \cite{jump} modeling language. \textcolor{black}{For CC OPF, we used the cutting-plane algorithm instead of the second-order cone reformulation. For RCC OPF, we used the algorithm described in Section~\ref{sec:robustmodel}.} \textcolor{black}{As shown in Figure~\ref{fig:cpu_time}, the} average wall-clock time for an RCC OPF instance was $\sim$8 seconds with an increase to $\sim$20 seconds for instances with $\Gamma$ = 1.0.  Such an increase is expected because of the larger number of cutting plane iterations required for constrained problems.

\textcolor{black}{This case study demonstrates that an instance can be solved within seconds on a system with 2209 buses, which is comparable to real-life power systems. Therefore, the proposed formulation is likely to be compatible with requirements of existing commercially available short-term planning tools.}

\subsection{Remarks}
\textcolor{black}{In the presented case study we made several assumptions which are specific to the BPA test system. As a result of these assumptions, the proposed methodology may lead to different cost performance and technical results, if applied to other systems. However,  we note that the BPA system takes advantage of highly flexible and low cost hydro generation, which are capable of accommodating variability and uncertainty of wind power generation at a relatively low cost. In other power systems, which predominantly feature fossil-fired and thus are less flexible and less cost effective in accommodating fluctuations of wind power generation, the proposed methodology is likely to be of greater value. We refer interested readers to \cite{ma_2013, ruth_2015}, which extensively discuss the impact of generation mix's flexibility and production cost on power system operations.}

\section{Conclusions}\label{sec:conclusions}

Based on \cite{ccopf-sirev}, we have developed and implemented algorithms to compute a distributionally-robust chance constrained optimal power flow (RCC OPF) that accounts for uncertainty in the parameters of statistical models that describe the deviations of wind (or other intermittent) generation from its forecast.

We have demonstrated the scalability of the RCC OPF by performing a seasonal case study on a modification of the BPA system. In this setting of vertically integrated grid operations, the case study shows that, compared to deterministic or even chance constrained OPF (CC OPF), the RCC OPF distributes both generation and regulation in a manner that can result in both cost savings and better technical performance; including fewer violations of transmission line limits, generator ramping limits, and smaller Area Control Error values.

The work in this manuscript points to several areas for potential future work, including:
\begin{itemize}
\item similar case studies should be performed on power systems that are dominated by fossil generation instead of hydro generation
\item the RCC OPF formulation should be extended to include the effects of reactive power on nodal voltage magnitudes, transmission line limits, and generator limits
\item to better model generator ramping constraints, the RCC OPF should be modified to a time-extended or look-ahead formulation consistent with the operation of modernized power systems
\item \textcolor{black}{the time-extended robust chance constraints should be adapted to day-ahead planning tools, such as UC and security-constrained UC. Further algorithmic developments may be needed to tractably solve such formulations.} 
\item the current formulation should be extended to market-based operations to incorporate the cost of procuring frequency regulation capacity.
\end{itemize}

\section*{Acknowledgement}
The work at LANL was funded by the Advanced Grid Modeling Program in the Office of Electricity in the US Department of
Energy and was carried out under the auspices of the National Nuclear Security Administration of the U.S. Department of
Energy at Los Alamos National Laboratory under Contract No. DE-AC52-06NA25396. M. Lubin was supported by the DOE Computational Science Graduate Fellowship, which is provided under grant number DE-FG02-97ER25308. Y. Dvorkin was supported in part by the Clean Energy Institute Student Training \& Exploration Grant.

\ifCLASSOPTIONcaptionsoff
  \newpage
\fi


\begin{IEEEbiographynophoto}{Miles Lubin} received his B.S. in Applied Mathematics and M.S. in Statistics from the University of Chicago in 2011. From 2011 to 2012 he was a predoctoral researcher at Argonne National Laboratory near Chicago, IL. He is currently a Ph.D. candidate in Operations Research at the Massachusetts Institute of Technology. Miles is a Department of Energy Computational Science Graduate Fellow and visitor at the Center for Nonlinear Studies at Los Alamos National Laboratory in 2014 and 2015. His research interests include large-scale mathematical optimization with application to power systems.
\end{IEEEbiographynophoto}

\begin{IEEEbiographynophoto}{Yury Dvorkin} (S'11) received the B.S.E.E degree with the highest honors at Moscow Power Engineering Institute (Technical University), Moscow, Russia, in 2011. He is currently pursuing the Ph.D. degree in electrical engineering at the University of Washington, Seattle, WA, USA.

Previously, Yury was a graduate intern at the Center for Nonlinear Studies at the Los Alamos National Laboratory, Los Alamos, NM, USA. He is a recipient of the Clean Energy Institute Graduate Fellowship (2013-2014) and the Clean Energy Institute Student Training \& Exploration Grant (2014-2015). His research interests include short- and long-term planning in power systems with renewable generation and power system economics.

\end{IEEEbiographynophoto}

\begin{IEEEbiographynophoto}{Scott Backhaus} received the Ph.D. degree in physics from the University of
California at Berkeley, Berkeley, CA, USA, in 1997 in the area of experimental
macroscopic quantum behavior of superfluid He-3 and He-4.

In 1998, he came to Los Alamos National Laboratory (LANL), Los Alamos, NM, USA, and was Director Funded Postdoctoral Researcher from 1998 to 2000, a Reines Postdoctoral Fellow from 2001 to 2003, and a Technical Staff Member from 2003 to now. While at LANL, he performed experimental and theoretical research in the area of thermoacoustic energy conversion.
Recently, his attention has shifted to other energy-related topics including the fundamental science of geologic carbon sequestration and grid integration of renewable generation.

Dr. Backhaus received an R\&D 100 award in 1999 and Technology Reviews Top 100 Innovators Under 35 (award in 2003) when he was with LANL. 
\end{IEEEbiographynophoto}

\vfill 
\end{document}